\documentclass[11pt]{article}
\usepackage[leqno]{amsmath}
\usepackage{amssymb}
\usepackage{amsthm}
\usepackage{amscd}

\newtheorem{thm}{Theorem}[section]
\newtheorem{cor}[thm]{Corollary}
\newtheorem{lem}[thm]{Lemma}
\newtheorem{conj}[thm]{Conjecture}

\newcommand{\set}[1]{\left\{#1\right\}}
\newcommand{\abs}[1]{\left\vert#1\right\vert}
\newcommand{\Irr}{\mathrm{Irr}}

\begin{document}

\title{Orders of elements and zeros and heights of characters in a finite group\footnote{2000 Mathematics Subject Classification 20C15, 20C20, 20F16.}}
\author{Tom Wilde\footnote{Correspondence: Tom Wilde, 2 Amner Road London SW116AA. Email tom@beech84.fsnet.co.uk.}}
\maketitle

\section{Introduction}
If $G$ is a finite group and $g\in G$ and $\chi\in\Irr(G)$ are such that $\chi(g)\neq 0,$ then by a fundamental result on
projective irreducible characters (e.g. \cite[Theorem 8.17]{Isaacs}), the squarefree part of the order of $g$ divides
$\abs{G}/\chi(1).$ We conjecture that the same holds for the order of $g$ itself:

\begin{conj}\label{conj1}
Let $G$ be a finite group. Let $\chi\in\Irr(G)$ and $g\in G$ and suppose that $\chi(g)\neq 0.$ Then the order of $g$ divides
$\abs{G}/\chi(1).$
\end{conj}

It is convenient to say that Conjecture \ref{conj1} holds \emph{at p} for a given prime $p,$ if we are able to prove that
$\chi(g)\neq 0$ implies that the $p$-part of the order of $g$ divides $\abs{G}/\chi(1).$ For example, if $\chi$ has height zero
in its $p$-block, then Conjecture \ref{conj1} holds at $p.$ To see this, let $D$ be a defect group of the $p$-block containing
$\chi.$ If $\chi(g)\neq 0$ then by e.g. \cite[Corollary 15.49]{Isaacs}, $g_p$ is conjugate to an element of $D.$ Hence the order
of $g_p$ divides $\abs{D}.$ But by definition, $\chi$ having height zero is equivalent to $\abs{D}=(\abs{G}/\chi(1))_p.$
Therefore $g_p$ divides $\abs{G}/\chi(1),$ as was to be shown. This example suggests a connection between Conjecture \ref{conj1}
and character heights; in fact we are able to exploit this to obtain some results on heights (see Section \ref{Heights}).

We briefly discuss evidence for Conjecture \ref{conj1}. It holds if $G$ is solvable, as we show in Section \ref{Solvable}. It
also holds for the symmetric groups, as a consequence of the hook length formula and Murnaghan-Nakayama rule \cite[\S4.12 and
\S4.45]{FultonHarris}. It is consistent with a theorem of Feit \cite[Theorem 1]{Feit}, which states that the values of $\chi$
are contained in the field of $\abs{G}/\chi(1)^{th}$ roots of unity. Empirically, we have verified it for all groups in the
GAP(2005) character table library \cite{GAP}; the library includes the ATLAS groups.\\

In this note, we present some results on Conjecture \ref{conj1}, together with related results and questions. In Section
\ref{Main} we prove a partial result for any finite group and in Section \ref{Heights} we give an application to character
heights. In Section \ref{Solvable} we prove Conjecture \ref{conj1} for solvable groups. In Section \ref{Centralchar} we
conjecture a congruence for central character values which implies Conjecture \ref{conj1} and is perhaps of some independent
interest. We prove this second conjecture in some special cases, including for solvable groups, and give some consequences.

\section{Results for arbitrary $G$}\label{Main}

In this section, we prove the following partial result for arbitrary finite groups $G:$

\begin{thm}\label{main}
Let $G$ be a finite group and suppose $g\in G$ and $\chi\in\Irr(G)$ are such that $\chi(g)\neq 0.$ Let $n$ be the order of $g,$
and let $n_0=\prod_{p\mid n}p$ be the squarefree part of $n.$ Then
$$\mathrm{(i)}\text{ }nn_0\text{ divides }(\frac{\abs{G}}{\chi(1)})^2\text{ and }\mathrm{(ii)}\text{ }n^3\text{ divides }\frac{\abs{G}^3}{\chi(1)^2}.$$
\end{thm}

Theorem \ref{main}$\mathrm{(i)}$ generalizes the classical fact that $n_0$ divides $\abs{G}/\chi(1)$ mentioned at the outset.
Theorem \ref{main} implies Conjecture \ref{conj1} when $\abs{G}$ is not divisible by any fifth power. More precisely:

\begin{cor}\label{cubefree}
In the situation of Theorem \ref{main}, if $p$ is a prime and $n_p\nmid\abs{G}/\chi(1),$ then $p^3$ divides $n,$ $p^3$ divides
$\chi(1),$ and $p^5$ divides $\abs{G}.$
\end{cor}

\begin{proof}[Proof of Corollary] The first two statements are immediate from Theorem \ref{main}$\mathrm{(i),(ii)}$ respectively.
Hence $p^{15}$ divides $n^3\chi(1)^2,$ so $p^5$ divides $\abs{G}$ by a second application of Theorem \ref{main}$\mathrm{(ii)}$.
\end{proof}

Our proof of Theorem \ref{main} is along similar lines to the proof of \cite[Theorem 8.17]{Isaacs}. We use the following
notation: Two elements of $G$ are \emph{rationally conjugate} if they generate conjugate cyclic subgroups of $G.$ If $x$ is a
$p$-element of $G,$ let $S_p(x)$ be the $p$-section of $x,$ that is the set of elements of $G$ whose $p$-part is conjugate to
$x,$ and write $T_p(x)$ for the \emph{rational} $p$-section of $x,$ that is the set of elements of $G$ whose $p$-part is
rationally conjugate to $x$. Let $R(G)$ be the ring of virtual characters of $G$ and let $A$ be the ring of algebraic integers.
For a subset $X\subseteq G,$ we write $\delta_X$ for the indicator function of $X,$ defined by $\delta_X(g)=1$ if $g\in X$ and
$\delta_X(g)=0$ otherwise. The first two lemmas are standard results.

\begin{lem}\label{1}
Let $x$ be a $p$-element of the finite group $G.$ If $S_p(x)$ and $T_p(x)$ are respectively the $p$-section and rational
$p$-section of $x,$ then $\abs{C_G(x)}_p\delta_{S_p(x)}\in A\otimes R(G)$ and $\abs{C_G(x)}_p\delta_{T_p(x)}\in R(G).$
\end{lem}

\begin{proof}
By Brauer's characterization of characters, it suffices to show that $\abs{C_G(x)}_p\delta_{S_p(x)}\vert_H\in A\otimes R(H)$ for
each nilpotent subgroup $H\subseteq G.$ But this is a linear combination of $\abs{C_H(x^\prime)}_p\delta_{S_p(x^\prime)}$ where
$x^\prime$ runs over a set of representatives for the $H$-conjugacy classes contained in $x^G\cap H.$ Hence we may assume $G=H$
is nilpotent. Writing $G=P\times Q$ where $P$ is the Sylow $p$-subgroup of $G,$ then $\delta_{S_p(x)}=\delta_{x^P}\times 1_Q,$
and the result follows as $C_P(x)\delta_{x^P}\in A\otimes R(P).$ The second assertion follows because $T_p(x)$ is a disjoint
union of $p$-sections, and if $\chi\in\Irr(G)$ then $[\chi,\delta_{T_p(x)}]$ is a rational number.
\end{proof}

\begin{lem}\label{2}
If $u$ is an integer of the cyclotomic field $E=\mathbb{Q}(\zeta_n)$ and $n_0$ is the squarefree part of $n,$ then
$tr_{E/\mathbb{Q}}(u)$ is divisible by $n/n_0.$
\end{lem}

\begin{proof}
Since the ring of integers of $E$ is $\mathbb{Z}[\zeta_n],$ we may assume $u=\zeta_m$ is a primitive $m^{th}$ root of unity for
some $m$ dividing $n$. Then
$$tr_{E/\mathbb{Q}}(u)=\sum_{0\leq r<n,(r,n)=1}\zeta_m^r=\frac{\varphi(n)}{\varphi(m)}\mu(m)$$
where $\varphi$ is Euler's function and $\mu$ is the M\"{o}bius function. However $\mu(m)$ is zero unless $m$ is square free, in
which case the right hand side is divisible by $\varphi(n)/\varphi(n_0)=n/n_0.$
\end{proof}

If $H$ is a subgroup of $G$ and $g\in G$ is an element of order $n,$ it is convenient to write $Aut_H(g)$ for the subgroup of
$(\mathbb{Z}/n\mathbb{Z})^*$ consisting of automorphisms induced on $\langle g\rangle$ by $N_H\langle g\rangle.$ Thus
$$Aut_H(g)=\set{r\in\mathbb{Z}/n\mathbb{Z},g^r\in g^H}$$Clearly $Aut_H(g)$ is isomorphic to $\frac{N_H\langle g\rangle}{C_H(g)}.$

\begin{lem}\label{3}
Suppose $H$ is a subgroup of $G.$ Let $g\in H$ have order $n,$ and let $X_g\subseteq H$ be the rational conjugacy class of $g$
in $H.$ If $\chi\in\Irr(G)$ and $\varphi\in\Irr(H)$ then $$\frac{1}{\abs{H}}\sum_{h\in
X_g}\chi(h)\bar{\chi}(h)\bar{\varphi}(h)\in\frac{\abs{C_G(g)}n/n_0}{(\abs{G}/\chi(1))^2\abs{Aut_H(g)}}\mathbb{Z}.$$
\end{lem}

\begin{proof}
The rational conjugacy class $X_g$ is the union of equivalence classes under the relation $h_1\sim h_2$ if $h_1$ and $h_2$
generate the same cyclic subgroup of $H.$ The number of these classes is $\abs{H:N_H\langle g\rangle}$ and for any one such
equivalence class $Y\subseteq X_g,$ it follows from the definition of the trace that $$\sum_{h\in
Y}\chi(h)\bar{\chi}(h)\bar{\varphi}(h)=tr_{\mathbb{Q}(\zeta_n)/\mathbb{Q}}(\chi(g)\bar{\chi}(g)\bar{\varphi}(g)).$$ Therefore
\begin{equation*}\begin{split}\frac{1}{\abs{H}}\sum_{h\in X_g}\chi(h)\bar{\chi}(h)\bar{\varphi}(h)
&=\frac{1}{\abs{H}}\frac{\abs{H}}{\abs{N_H\langle g\rangle}}tr_{\mathbb{Q}(\zeta_n)/\mathbb{Q}}(\chi(g)\bar{\chi}(g)\bar{\varphi}(g))\\
&=\frac{1}{\abs{N_H\langle g\rangle}}\frac{\abs{C_G(g)}^2}{(\abs{G}/\chi(1))^2}tr_{\mathbb{Q}(\zeta_n)/\mathbb{Q}}(\omega_\chi(g)\bar{\omega}
_\chi(g)\bar{\varphi}(g))\\
&=\frac{\abs{C_G(g)}}{\abs{C_H(g)}}\frac{\abs{C_G(g)}}{(\abs{G}/\chi(1))^2\abs{Aut_H(g)}}
tr_{\mathbb{Q}(\zeta_n)/\mathbb{Q}}(\omega_\chi(g)\bar{\omega}_\chi(g)\bar{\varphi}(g))\\
&\in\frac{\abs{C_G(g)}n/n_0}{(\abs{G}/\chi(1))^2\abs{Aut_H(g)}}\mathbb{Z},
\end{split}\end{equation*}
where $\omega_\chi(g)=\frac{h_g\chi(g)}{\chi(1)}$ is the value of the central character associated with $\chi,$ and the last
step follows by Lemma \ref{2}, because $\omega_\chi(g)\bar{\omega}_\chi(g)\bar{\varphi}(g)$ is an algebraic integer.
\end{proof}

If $g\in G,$ let $o(g)$ denote the order of $g$ and let $\nu_p(n)$ be the power of $p$ which divides the integer $n.$ Theorem
\ref{main} will follow easily once we have proved the next result.

\begin{thm}\label{AutIneq}
Let $\chi\in\Irr(G).$ Suppose $g\in G$ has $\chi(g)\neq 0.$ Then for any prime $p,$
$$2\nu_p(\frac{\abs{G}}{\chi(1)})+\nu_p(\abs{Aut_G(g_p)})\geq 2\nu_p(o(g)).$$
\end{thm}

\begin{proof}
Let $x=g_p$ and define a class function $\theta$ on $G$ by $\theta=\chi\bar{\chi}\delta_{T_p(x)}.$ By hypothesis, $\chi(g)\neq
0,$ so $\theta$ is not identically zero. We may clearly suppose $x\neq 1,$ so that $1\notin T_p(x).$ Hence
$$0<[\theta,1_G]=\frac{1}{\abs{G}}\sum_{h\in T_p(x)}\abs{\chi(h)}^2<[\chi,\chi]=1,$$ so $\theta$ is not a virtual
character of $G.$ By Lemma \ref{1}, $\abs{G}_p\theta\in R(G),$ so we must have $m\theta\notin R(G)$ whenever $m$ is prime to
$p.$ By Brauer's characterization of characters (local form), there exists a $p$-elementary subgroup $H\subseteq G$ and a
character $\varphi\in\Irr(H)$ such that $\abs{H}_{p^\prime}[\theta,\varphi]_H\notin\mathbb{Z}.$ Then $T_p(x)\cap H$ cannot be
empty, and is a union of rational conjugacy classes of $H,$ and for one of these, say $X_{g^\prime}\subseteq H,$
$\frac{1}{\abs{H}_p}\sum_{h\in X_{g^\prime}}\chi(h)\bar{\chi}(h)\bar{\varphi}(h)\notin\mathbb{Z}.$ By Lemma \ref{3}, we must
have
$$2\nu_p(\frac{\abs{G}}{\chi(1)})+\nu_p(\abs{Aut_H(g^\prime)})>\nu_p(\abs{C_G(g^\prime)})+\nu_p(o(g^\prime))-1.$$ Since $H$ is
$p$-elementary, $\abs{Aut_H(g^\prime)}=\abs{Aut_H(g^\prime_p)},$ and since $g^\prime_p$ and $g_p$ are rationally conjugate, we
have $\nu_p(o(g^\prime))=\nu_p(o(g)),$ $\nu_p(\abs{C_G(g^\prime)})\geq\nu_p(o(g))$ and $\abs{Aut_H(g^\prime_p)}$ divides
$\abs{Aut_G(g_p)}.$ The desired inequality follows.
\end{proof}

\begin{proof}[Proof of Theorem \ref{main}]
If $g$ has order $n$ and $k=\nu_p(n)\geq 1$ then $Aut_G(g_p)$ is a subgroup of $(\mathbb{Z}/p^k\mathbb{Z})^\star,$ which has
order $p^{k-1}(p-1).$ Therefore $\nu_p(\abs{Aut_G(g_p)})\leq\nu_p(n)-1.$ Hence by Theorem \ref{AutIneq}, if $\chi(g)\neq 0$ then
$2\nu_p(\abs{G}/\chi(1))\geq\nu_p(n)+1.$ This holds for each prime dividing $n,$ so $nn_0$ divides $(\abs{G}/\chi(1))^2,$ which
is Theorem \ref{main}$\mathrm{(i)}$. On the other hand, $\nu_p(\abs{Aut_G(g_p)})\leq\nu_p(\abs{G})-\nu_p(n).$ Hence by Theorem
\ref{AutIneq}, if $\chi(g)\neq 0$ then $2\nu_p(\abs{G}/\chi(1))+\nu_p(\abs{G})\geq 3\nu_p(n),$ which similarly implies Theorem
\ref{main}$\mathrm{(ii)}$.
\end{proof}

In the proof of Theorem \ref{main}, the term $\nu_p(\abs{Aut_G(g_p)})$ from Theorem \ref{AutIneq} is the obstruction to a proof
of Conjecture \ref{conj1}. If for a given prime $p,$ the Sylow $p$-subgroups of $G$ are abelian, this term vanishes, and the
proof of Theorem \ref{main} establishes that $n_p$ divides $\abs{G}/\chi(1),$ so Conjecture \ref{conj1} holds at $p.$ Using the
following theorem of Willems, we can obtain the same conclusion if a defect group of the $p$-block containing $\chi$ is abelian.

\begin{thm}\label{Willems}{\rm(\cite[Corollary 3.5]{Willems})}
Let $\chi\in\Irr(B)$ where $B$ is a $p$-block of the finite group $G.$ Let $D$ be a defect group of $B.$ Then there exist
$a_i\in\mathbb{Z}$ and elementary subgroups $H_i\subseteq G$ such that $$\chi=\sum a_i\lambda_i^G,$$ where $\lambda_i$ is a
linear character of $H_i$ and the Sylow $p$-subgroup of each $H_i$ is contained in $D.$
\end{thm}

\begin{thm} \label{sub1}
In the situation of Theorem \ref{main}, let $p$ be a prime dividing $n$ and let $B$ be the block of $G$ in characteristic $p$
which contains $\chi.$ If the defect groups of $B$ are abelian, then $n_p$ divides $\abs{G}/\chi(1).$ In other words, Conjecture
\ref{conj1} holds for $\chi$ at $p.$
\end{thm}

\begin{proof}
Let $D$ be a defect group of $B$ and let $\chi=\sum a_i\lambda_i^G,$ be the sum in Theorem \ref{Willems}. Let
$\theta=\chi\bar{\chi}\delta_{T_p(x)}$ be the class function defined in the proof of Theorem \ref{AutIneq}, and also write
$\theta_1=\chi\delta_{T_p(x)}.$ As in the proof of Theorem \ref{AutIneq}, $m\theta\notin R(G)$ for any $m$ prime to $p.$ Since
$$\theta=\bar{\chi}\theta_1=\sum a_i(\bar{\lambda}_i\theta_1\vert_{H_i})^G,$$ it follows that for some $i$ there exists a
$p$-elementary subgroup $H$ with $H\subseteq H_i$ and a character $\varphi\in\Irr(H)$ such that
$\abs{H}_{p^\prime}[\theta_1,\varphi]_H\notin\mathbb{Z}.$ Again $T_p(x)\cap H$ cannot be empty, and is a union of rational
conjugacy classes of $H,$ and for one of these, say $X_{g^\prime}\subseteq H,$ $\frac{1}{\abs{H_p}}\sum_{h\in
X_{g^\prime}}\chi(h)\bar{\varphi}(h)\notin\mathbb{Z}.$ This is just as in the proof of Theorem \ref{AutIneq}, except that we
have $\chi$ in place of $\chi\bar\chi.$ The proof of Lemma \ref{2} with $\chi$ in place of $\chi\bar{\chi}$ shows that
$$\frac{1}{\abs{H}_p}\sum_{h\in
X_{g^\prime}}\chi(h)\bar{\varphi}(h)\in\frac{n^\prime/n^\prime_0}{(\abs{G}/\chi(1))\abs{Aut_H({g^\prime})}}\mathbb{Z},$$ where
$n^\prime=o(g^\prime).$ But $H$ is $p$-elementary, and its Sylow $p$-subgroup is contained in $D,$ so in fact $H$ is abelian and
$Aut_H(g^\prime)=1.$ We conclude that $\nu_p(\abs{G}/\chi(1))\geq\nu_p(n^\prime)=\nu_p(n),$ as required.
\end{proof}

If $D$ is abelian, then according to Brauer's height zero conjecture, $\chi$ in Theorem \ref{sub1} should have height zero, from
which it would follow that Conjecture \ref{conj1} holds at $p,$ as remarked in the introduction. The above proof does not rely
on the height zero conjecture; moreover, the fact that $D$ is abelian was only used to ensure $Aut_H(g^\prime)=1.$ The proof
therefore goes through unchanged provided that in $D,$ each cyclic subgroup is in the center of its normalizer. This holds, for
example, if $D$ is the direct product of an abelian group and a group of exponent $p.$

\section{An application to heights of characters}\label{Heights}

The results in Section \ref{Main} can be used to give some control over heights of characters. If $B$ is a $p$-block of $G,$ we
write $ht_p(\chi)$ for the height of $\chi$ in $B,$ so that if $\abs{G}_p=p^a$ and $B$ has defect $d$ then
$\chi(1)_p=p^{a-d+ht_p(\chi)}.$ For any group $G,$ write $e(G)$ for the exponent of $G.$

\begin{thm}\label{Ht1}
Let $B$ be a $p$-block of the finite group $G$ and let $D$ be a defect group of $B.$ Let $\abs{G}_p=p^a$ and $\abs{D}=p^d.$ Let
$\chi\in\Irr(B).$ \\\\ $\mathrm{(i)}$ If Conjecture \ref{conj1} holds for $G$ at $p,$ then $ht_p(\chi)\leq d-\nu_p(e(Z(D)).$ In
particular, this holds if $D$ is abelian or $G$ is solvable.
\\\\$\mathrm{(ii)}$ Without assuming Conjecture \ref{conj1}, we have $ht_p(\chi)\leq\frac{a+d}{2}-\nu_p(e(Z(D)).$ If $a-d\leq 1$
then $\mathrm{(i)}$ follows without Conjecture \ref{conj1}. In particular, $(\rm{i})$ holds if $B$ is the principal block of any
finite group.
\end{thm}

\begin{proof}
Let $D$ be a defect group of $B,$ and let $y$ be such that $y^G$ is a defect class of $B$ and $D$ is a Sylow $p$-subgroup of
$C_G(y).$ Let $x\in Z(D)$ and let $g=xy.$ We claim that $\chi(g)\neq 0.$ This follows from \cite[Lemma 5.15(b)]{Navarro}, but we
sketch the short argument. First, $D$ is a Sylow $p$-subgroup of $C_G(g),$ so if $\psi\in\Irr(B)$ is a character of height zero,
then $h_g/\psi(1)$ is an integer prime to $p.$ Also if $\mathfrak{p}$ is a prime lying over $p$ in a splitting field for $G,$
then $\psi(g)=\psi(y)\neq 0\mod\mathfrak{p},$ since $y$ is a defect class of $B.$ Hence
$$\omega_\psi(g)=\frac{h_{g}\psi(g)}{\psi(1)}\neq 0\mod\mathfrak{p}.$$ Since $\chi$ and $\psi$ lie in the same block, also
$\omega_\chi(g)\neq 0\mod\mathfrak{p},$ so certainly $\chi(g)\neq 0$ as required. Now assuming Conjecture \ref{conj1} at $p,$ we
conclude that $o(x)$ divides $(\abs{G}/\chi(1))_p=p^{d-ht_p(\chi)}.$ Since $x\in Z(D)$ was arbitrary, $(\mathrm{i})$ follows.
Conjecture \ref{conj1} is proved at Theorem \ref{sub1} for $D$ abelian, and at Corollary \ref{ConjSolvable} below for $G$
solvable.

Next, we do not assume Conjecture \ref{conj1} and appeal instead to Theorem \ref{AutIneq}. Since $x\in Z(D),$ we have
$\nu_p(\abs{Aut_G(x)})\leq a-d,$ so by Theorem \ref{AutIneq}, $$2(d-ht_p(\chi))+a-d\geq 2\nu_p(o(x)).$$ As $x\in Z(D)$ was
arbitrary, this gives $(\mathrm{ii}).$ If $a-d\leq 1$ then $\frac{a+d}{2}\leq d+\frac{1}{2},$ so $(\mathrm{ii})$ implies
$(\mathrm{i})$ in this case.
\end{proof}

The solvable case of Theorem \ref{Ht1}$(\mathrm{i})$ is really redundant, since the better inequality
$ht_p(\chi)\leq\nu_p(\abs{D:Z(D)})$ is true for $p$-solvable groups by a well known result of Fong (\cite[Theorem 3C]{Fong}, or
see \cite[Theorem 10.21]{Navarro}). However, if $G$ is not solvable, then Theorem \ref{Ht1} seems in some cases to give a better
bound than known results. For example, the bound in \cite[Corollary $9.11$]{Feit2} is $ht_p(\chi)<
d-\nu_p(\abs{Z(D)})+\frac{1}{2}r(r+1)$ where $r=\nu_p(\abs{D:\Phi(D)}).$ If say $D=C_{p^e}\times (C_p)^{r-1}$ then this gives
$ht_p(\chi)<\frac{1}{2}r(r+1),$ while Theorem \ref{Ht1} gives $ht_p(\chi)< r.$

\section{Solvable groups}\label{Solvable}

Conjecture \ref{conj1} is true for solvable groups. More generally, if $G$ is $p$-solvable then Conjecture \ref{conj1} holds at
$p$ for $G,$ in the sense mentioned in the introduction. This will follow easily from:

\begin{thm}\label{PrimitivepSolvable}
Let $G$ be a finite $p$-solvable group. Suppose $\chi\in\Irr(G)$ is primitive. If $Q$ is an abelian $p$-subgroup of $G$ then
$\chi(1)$ divides $\abs{G:Q}.$
\end{thm}

\begin{proof}
We use some theory of special characters of $p$-solvable groups (see for example \cite[Chapter VI]{ManzWolf}). Namely, as $\chi$
is primitive, there is a factorization $\chi=\chi_p\chi_{p^\prime}$ where $\chi_p$ and $\chi_{p^\prime}$ are $p$-special and
$p^\prime$-special respectively. Let $P$ be a Sylow $p$-subgroup of $G$ containing $Q.$ We must show that $\chi_p(1)=\chi(1)_p$
divides $\abs{P:Q}.$ But $\chi_p$ remains irreducible on restriction to $P,$ so since $Q$ is abelian, $\chi_p(1)\leq\abs{P:Q},$
and the result follows as $\chi_p(1)$ and $\abs{P:Q}$ are powers of $p.$
\end{proof}

Theorem \ref{PrimitivepSolvable} does not hold for all groups. For example, the alternating group $A_5$ has a primitive
character of degree $4.$

\begin{cor}\label{PrimitiveSolvable}
In the situation of Theorem \ref{PrimitivepSolvable}, if $G$ is solvable and $A\subseteq G$ has abelian Sylow subgroups, then
$\chi(1)$ divides $\abs{G:A}$
\end{cor}

\begin{proof}
Immediate from Theorem \ref {PrimitivepSolvable}.
\end{proof}

Corollary \ref{PrimitiveSolvable} appears as \cite[Theorem 5.15b]{Trento} under the (unnecessary) hypothesis that $A$ is normal
in $G.$ Our proof above is essentially the same as the proof given there.

\begin{cor}\label{ConjSolvable}
Conjecture \ref{conj1} holds for solvable groups.
\end{cor}

\begin{proof}
Let $\chi\in\Irr(G)$ where $G$ is solvable, and let $\varphi\in\Irr(H)$ be a primitive character of $H\subseteq G$ such that
$\chi=\varphi^G.$ If $g\in G$ and $\chi(g)\neq 0$ then a conjugate of $g$ lies in $H.$ Then Corollary \ref{PrimitiveSolvable}
implies that the order of $g$ divides $\abs{H}/\varphi(1)=\abs{G}/\chi(1).$
\end{proof}

It would be interesting if Theorem \ref{PrimitivepSolvable} could be used along the lines of Theorem \ref{Ht1} to prove Fong's
theorem that $ht_p(\chi)\leq\nu_p(\abs{D:Z(D)})$ for a character in a $p$-block of a $p$-solvable group with defect group $D.$
We show this can be done in the easiest case, when $B$ is the principal block of $G.$ We fix a prime $p$ and write $G^0$ for the
set of $p$-regular elements of $G$ and $[\chi,\psi]^0$ for $\frac{1}{\abs{G}}\sum_{g\in G^0}\chi(g^{-1})\psi(g).$

\begin{cor}
Let $G$ be a finite $p$-solvable group and suppose $\chi\in\Irr(G)$ belongs to the principal $p$-block of $G.$ Let $P$ be a
Sylow $p$-subgroup of $G.$ Then $ht_p(\chi)\leq\nu_p(\abs{P:Z(P)}).$
\end{cor}

\begin{proof}
As is well known, $\chi$ is lifted from a character in the principal block of $G/\mathrm{O}_{p^\prime}(G),$ so by induction we
assume $\mathrm{O}_{p^\prime}(G)=1.$ We have to prove that $\chi(1)_p$ divides $\abs{G:Z(D)}.$ This follows from Theorem
\ref{PrimitivepSolvable} if $\chi$ is primitive, so we may assume that $\chi=\varphi^G$ for some maximal subgroup $H$ of $G$ and
$\varphi\in\Irr(H).$ Since $\chi$ is in the principal block, $[\chi,1_G]^0_G\neq 0,$ so $[\varphi,1_H]^0_H\neq 0$ by Frobenius
reciprocity; therefore $\varphi$ belongs to the principal block of $H$ (or use the Third Main Theorem). By induction,
$\varphi(1)$ divides $\abs{H:Z(Q)}$ where $Q$ is a Sylow $p$-subgroup of $H.$ Hence, it suffices to show that $Z(P)\subseteq H.$
Suppose the contrary and let $x\in Z(P)\backslash H.$ Since $\mathrm{O}_{p^\prime}(G)=1,$ the Hall-Higman Lemma implies $x\in
Z(\mathrm{O}_p(G)).$ In particular, $\mathrm{O}_p(G)\nsubseteq H,$ and as $H$ is maximal we have $G=H\mathrm{O}_p(G).$ Then
since $x\in Z(\mathrm{O}_p(G)),$ we have $x^G=x^H,$ so $x^G$ is disjoint from $H$ and we conclude that $\chi(x)=0.$ As $\chi$ is
in the principal block, $0=\omega_\chi(x)=\abs{h_x}\mod\mathfrak{p},$ where $\mathfrak{p}$ is a prime lying over $p$ in a
splitting field for $G.$ Hence $p\mid h_x,$ contradicting the fact that $x\in Z(P).$ This contradiction shows that
$Z(P)\subseteq H$ and the proof is complete.
\end{proof}

We mention a relative of Corollary \ref{PrimitiveSolvable} in which, using projective representations, we replace the character
degree $\chi(1)$ by the ramification index of $\chi$ with respect to a normal subgroup.

\begin{thm}\label{PrimitiveSolvable2}
Let $G$ be a finite group and suppose $\chi\in\Irr(G)$ is primitive. Let $N\vartriangleleft G$ be a normal subgroup of $G$ and
suppose $G/N$ is solvable. Let $\varphi\in\Irr(N)$ be the character of $N$ lying under $\chi.$ Let $e(G/N)$ be the exponent of
$G/N.$ Then $\chi(1)/\varphi(1)$ divides $\abs{G:N}/e(G/N).$
\end{thm}

\begin{proof}
Since $\chi$ is primitive, $\chi_N=f\varphi$ where $f$ is the ramification index. By the theory of projective representations,
there exists a central extension $$1\rightarrow Z\rightarrow\Gamma\xrightarrow{\pi} G\rightarrow 1$$ such that
$\pi^{-1}(N)=\tilde{N}\times Z$ for some normal subgroup $\tilde{N}\vartriangleleft \Gamma,$ and $\Gamma$ has irreducible
characters $\alpha,\beta\in\Irr(\Gamma)$ with $\beta(1)=f,$ $\tilde{N}\subseteq\ker\beta,$ and $Z\subseteq\ker(\alpha\beta)$ and
$\alpha\beta=\chi$ viewed as a character of $G.$ We claim that $\beta$ is primitive. For suppose not, so that
$\beta=\eta^\Gamma$ where $\Sigma\subset\Gamma$ is a proper subgroup of $\Gamma$ and $\eta\in\Irr(\Sigma).$ Then
$\alpha\beta=\alpha\eta^\Gamma=(\alpha_\Sigma\eta)^\Gamma.$ Also $\ker(\alpha\beta)=\mathrm{core}_G(\ker(\alpha_\Sigma\eta))$ by
\cite[Lemma 5.11]{Isaacs}, so $Z\subseteq\Sigma$ and $\alpha_\Sigma\eta$ is the lift of some $\psi\in\Irr(H),$ where $H$ is the
image of $\Sigma$ in $G.$ But then $H$ is a proper subgroup of $G$ and $\chi=\psi^G,$ contrary to hypothesis. ($\alpha$ is
primitive for the same reason but we do not use this.)

Now the sequence above gives rise to the deflated sequence $$1\rightarrow
Z\rightarrow\Gamma/\tilde{N}\xrightarrow{\tilde{\pi}}G/N\rightarrow 1.$$ If $C\subseteq G/N$ is any cyclic subgroup then
$\tilde{\pi}^{-1}(C)$ is abelian. $\beta$ is the inflation of a character of $\Gamma/\tilde{N},$ which must be primitive, so by
Corollary \ref{PrimitiveSolvable}, $\beta(1)$ divides $\vert\Gamma/\tilde{N}:\tilde{\pi}^{-1}(C)\vert=\vert G/N:C\vert,$ and the
result follows.
\end{proof}

For solvable $G,$ Corollary \ref{PrimitiveSolvable} is an essentially stronger result than Conjecture \ref{conj1}. More is also
known in a slightly different direction. The following is the main result of \cite{Wilde}:

\begin{thm}\label{BrauerInd}{\rm(\cite{Wilde}, Main Theorem)}
Let $G$ be a finite solvable group and suppose $\chi\in\Irr(G).$ Then there exist subgroups $H_i\subseteq G$ such that $\chi(1)$
divides $\abs{G:H_i}$ for each $i,$ and a sum $$\chi=\sum_ia_i\lambda_i^G$$ where for each $i,$ $a_i$ is an integer and
$\lambda_i$ is a linear character of $H_i.$
\end{thm}

Conjecture \ref{conj1} follows from Theorem \ref{BrauerInd}, since if $g\in G$ and $\chi(g)\neq 0$ then a conjugate of $g$ lies
in some $H_i,$ and then the order of $g$ divides $\abs{H_i}$ which divides $\abs{G}/\chi(1).$ Unlike Theorem
\ref{PrimitiveSolvable}, it seems possible that Theorem \ref{BrauerInd} holds for all finite groups. We remark that, for any $G$
and $\chi\in\Irr(G),$ if $\chi$ has height zero in its $p$-block, then the sum in Theorem \ref{Willems} has $\chi(1)_p$ dividing
$\abs{G:H_i},$ so that in the obvious sense, Theorem \ref{BrauerInd} holds for $\chi$ at $p.$ \\

Lastly in this section, we remark that the conclusion of Theorem \ref{BrauerInd} also holds for any group $G$ when
$\abs{G}/\chi(1)$ is a prime power. In fact, then $\chi$ is monomial:

\begin{thm}
Let $\chi\in\Irr(G),$ where $G$ is a finite group, and suppose that $\abs{G}/\chi(1)=p^k$ where $p$ is prime. Then $\chi$ is
monomial.
\end{thm}

\begin{proof}
Let $P$ be a Sylow $p$-subgroup of $G$ and let $\abs{P}=p^m.$ Then $\chi(1)=p^{m-k}\abs{G}_{p^\prime}.$ Since
$p^{m-k+1}\nmid\chi(1),$ some irreducible component $\psi$ of $\chi\vert_P$ has degree $\psi(1)\leq p^{m-k}.$ But then
$\psi^G(1)\leq p^{m-k}\abs{G}_{p^\prime}\leq\chi(1),$ so since $[\chi,\psi^G]\neq 0,$ equality holds and $\chi=\psi^G.$ The
result follows since $\psi$ is monomial.
\end{proof}

\section{Central character values}\label{Centralchar}

In our proof of Theorem \ref{main}, the term $\nu_p(\abs{Aut_G\langle g_p\rangle})$ from Theorem \ref{AutIneq} prevents us from
establishing Conjecture \ref{conj1}. This term would be removable if in Lemma \ref{3}, $\omega_\chi(g)$ were divisible by
$\abs{Aut_H(g)}.$ However, it is not true in general that $\abs{Aut_G(g)}$ divides $\omega_\chi(g)$ for any $\chi\in\Irr(G)$ and
$g\in G.$ For example, let $\chi$ be the permutation character of degree $4$ of the symmetric group $S_5$ and let $g$ be a
$5$-cycle. Then $\omega_\chi(g)=-6$ is not divisible by $\abs{Aut_G(g)}=4.$ Solvable examples also exist.

However, there is some evidence that $\omega_\chi(g)$ is divisible by the order of a certain large subgroup $Aut^0_G(g)\subseteq
Aut_G(g).$ If $H$ is a subgroup of $G$ and $g\in G$ has order $n,$ then we define $Aut^0_H(g)\subseteq Aut_H(g)$ as follows. For
each prime $p$ dividing $n,$ set
$$Aut^0_{H,p}(g)=\set{m\in Aut_H(g);m=1 \mod pn_{p^\prime}}$$ if $p$ is odd or $p=2$ and $4\nmid n$, and
$$Aut^0_{H,2}(g)=\set{m\in Aut_H(g);m=1 \mod 4n_{2^\prime}}$$ if $p=2$ and $4\mid n.$ Finally, set
$$Aut^0_H(g)=\prod_{p\mid n}Aut^0_{H,p}(g).$$

If $p$ is odd then $Aut^0_{H,p}(g)$ is naturally isomorphic to the Sylow $p$-subgroup of $\frac{N_H\langle g\rangle\cap
C_H(g_{p^\prime})}{C_H(g)}.$ If $p=2$ and $4$ divides $n$ then $Aut^0_{H,p}(g)$ is naturally isomorphic to the subgroup of index
at most $2$ in this Sylow $p$-subgroup, consisting of elements that centralize the unique subgroup of order $4$ in $\langle
g\rangle.$ Finally, note that if $H\subseteq G$ is nilpotent and $g\in H$ then $Aut^0_H(g)\subseteq Aut_H(g)$ has index at most
$2.$

\begin{conj}\label{conj2}
Let $\chi\in\Irr(G)$ and for $g\in G,$ let $\omega_\chi(g)=\frac{h_g\chi(g)}{\chi(1)}$ be the central character value. Then
$\omega_\chi(g)=0\mod \vert Aut^0_G(g)\vert.$
\end{conj}

As with Conjecture \ref{conj1} it is convenient to have a notion of Conjecture \ref{conj2} holding \emph{at p} for a prime $p.$
Here, this means that for all $g\in G,$ $\omega_\chi(g)$ is divisible by $\vert Aut^0_{G,p}(g)\vert=\vert Aut^0_G(g)\vert_p.$

We will prove in Corollary \ref{conjsolvable} that Conjecture \ref{conj2} is true when $G$ is solvable, and is true at $p$
whenever $\chi$ has height zero in its $p$-block. We first show that Conjecture \ref{conj2} implies Conjecture \ref{1}.
Specifically:

\begin{thm}\label{1implies2}
Let $G$ be a finite group and suppose $\chi\in\Irr(G).$ If Conjecture \ref{conj2} holds for $\chi$ then Conjecture \ref{conj1}
also holds for $\chi.$
\end{thm}

\begin{proof}
Assume that $\chi\in\Irr(G)$ and $g\in G$ are such that $\chi(g)\neq 0,$ and let $x$ be the $p$-part of $g$ for a fixed prime
$p.$ Following the proof of Theorem \ref{AutIneq}, there is a $p$-elementary subgroup $H\subseteq G$ and a rational conjugacy
class $X_{g^\prime}\subseteq T_p(x)\cap H,$ such that $\frac{1}{\abs{H}_p}\sum_{h\in
X_{g^\prime}}\chi(h)\bar{\chi}(h)\bar{\varphi}(h)\notin\mathbb{Z}.$ Assuming Conjecture \ref{conj2}, we may strengthen Lemma
\ref{3} to
$$\frac{1}{\abs{H}}\sum_{h\in X_g}\chi(h)\bar{\chi}(h)\bar{\varphi}(h)\in\frac{\abs{C_G(g)}n/n_0\abs{Aut^0_G(g)}^2}{(\abs{G}/
\chi(1))^2\abs{Aut_H(g)}}\mathbb{Z}.$$ Hence as in the proof of Theorem \ref{AutIneq}, we now obtain:
$$2\nu_p(\frac{\abs{G}}{\chi(1)})+\nu_p(\abs{Aut_H(g^\prime)})-2\nu_p(\abs{Aut^0_G(g^\prime)})>\nu_p(\abs{C_G(g^\prime)})+\nu_p(o(g^\prime))-1.$$
Since $H$ is nilpotent, $Aut^0_H(g)\subseteq Aut_H(g)$ has index at most $2,$ so $\nu_p(\abs{Aut_H(g^\prime)})\leq
\nu_p(\abs{Aut^0_G(g^\prime)})+1$ (with the $1$ only needed when $p=2$) and
$$2\nu_p(\frac{\abs{G}}{\chi(1)})+1-\nu_p(\abs{Aut^0_G(g^\prime)})\geq \nu_p(\abs{C_G(g^\prime)})+\nu_p(o(g^\prime)).$$ Hence
certainly $$2\nu_p(\frac{\abs{G}}{\chi(1)})+1\geq \nu_p(\abs{C_G(g^\prime)})+\nu_p(o(g^\prime)),$$ and arguing as in the proof
of Theorem \ref{AutIneq}, we obtain $2\nu_p(\frac{\abs{G}}{\chi(1)})+1\geq 2\nu_p(o(g)).$ The $1$ may be discarded, and
Conjecture \ref{conj1} follows.
\end{proof}

In the penultimate equation of the above proof, we replaced $\nu_p(\abs{Aut^0_G(g^\prime)})$ with zero. We cannot make better
use of this term in general because we do not have control over the $p^\prime$-part of $g^\prime.$ If $G$ is rational, however,
then we have $\abs{Aut^0_G(g)}=n/2n_0$ if $n=o(g)$ is divisible by $4,$ and $\abs{Aut^0_G(g)}=n/n_0$ otherwise. Since the
$p$-parts of $g$ and $g^\prime$ are rationally conjugate, we know in this case that $\nu_p(\abs{Aut^0_G(g^\prime)})\geq
\nu_p(n)-2$ for $p=2,$ and $\geq \nu_p(n)-1$ for odd $p.$ The penultimate equation in the proof above is therefore improved to
$$2\nu_p(\frac{\abs{G}}{\chi(1)})+1-\nu_p(o(g))+2\geq 2\nu_p(o(g)),$$ where the $1$ and $2$ can be replaced with $0$ and $1$
respectively for $p$ odd. We have shown:

\begin{thm}\label{rationalgroup}
Let $G$ be a finite rational group, let $\chi\in\Irr(G)$ and suppose Conjecture \ref{conj2} holds for $G$ and $\chi.$ If $g\in
G$ has order $n$ and $\chi(g)\neq 0,$ then $n^3/n_0\text{ divides }4(\abs{G}/\chi(1))^2.$
\end{thm}

In particular, Theorem \ref{rationalgroup} will hold when $G$ is a rational solvable group, by Theorem \ref{conjsolvable}.
Another result depending on Conjecture \ref{conj2}, and so true for solvable groups $G$ is the following.

\begin{thm}
Let $G$ be a finite group. Suppose $\chi\in\Irr(G)$ is a faithful irreducible character with $\chi(1)$ a power of an odd prime
$p.$ Suppose further that Conjecture \ref{conj2} holds for $G$ and $\chi.$ If $P$ is a Sylow $p$-subgroup of $G$ and $g\in P$ is
such that $\langle g\rangle \vartriangleleft P,$ then either $g\in Z(G)$ or $\chi(g)=0.$
\end{thm}

\begin{proof}
Since $Aut_P(g)=P/C_P(g)$ and $C_P(g)$ is a Sylow subgroup of $C_G(g),$ the $p$-part of $h_g$ equals $\abs{Aut_P(g)}.$ Also
$Aut_P(g)\subseteq Aut^0_G(g),$ so Conjecture \ref{conj2} implies that $\chi(1)$ divides $\chi(g).$ Hence by the standard
argument of Burnside, either $\chi(g)=0$ or $g\in Z(G/\ker\chi)=Z(G).$
\end{proof}

We turn to the proof of Conjecture \ref{conj2} for solvable groups. The following well known lemma is a $p$-local form of Lemma
\ref{2}.

\begin{lem}\label{tracelemma}
Suppose $n=p^kn_{p^\prime}$ where $p$ is odd and $k\geq 2$ and suppose $E$ and $F$ are fields with
$\mathbb{Q}(\zeta_{pn_{p^\prime}})\subseteq F\subseteq E\subseteq\mathbb{Q}(\zeta_n).$ If $p=2,$ suppose additionally that
$i=\sqrt{-1}\in F.$ If $u\in E$ is an algebraic integer then $tr_{E/F}(u)$ is divisible by $\abs{E:F}.$
\end{lem}

\begin{proof}
First suppose $p$ is odd. Then $\mathbb{Q}(\zeta_n)/\mathbb{Q}(\zeta_{pn_{p^\prime}})$ is a cyclic Galois extension of prime
power degree having the unique chain of subfields
$$\mathbb{Q}(\zeta_{pn_{p^\prime}}=\zeta_n^{p^{k-1}})\subset...\subset\mathbb{Q}(\zeta_n^p)\subset\mathbb{Q}(\zeta_n).$$ $E$ is
one of these fields, so it is no loss to assume that $E=\mathbb{Q}(\zeta_n)$ and $F=\mathbb{Q}(\zeta_n^{p^r})$ for some $r$ with
$0\leq r\leq k-1.$ Then $\abs{E:F}=p^r,$ and the Galois group of $E/F$ is $(1+p^{k-r}n_{p^\prime}\mathbb{Z})/n\mathbb{Z}.$ Since
the ring of integers of $E$ is $\mathbb{Z}[\zeta_n],$ we may assume that $u=\zeta_n^m$ is a root of unity. Then
$$tr_{E/F}(u)=\sum_{s=0}^{p^r-1}(\zeta_n^m)^{1+p^{k-r}n_{p^\prime}s}=\zeta_n^m\sum_{s=0}^{p^r-1}\omega^s,$$ where
$\omega=(\zeta_n^m)^{p^{k-r}n_{p^\prime}}.$ If $\omega=1$ then the sum is divisible by $p^r=\abs{E:F};$ otherwise the sum is
$\zeta_n^m(\omega^{p^r}-1)/(\omega-1)=0.$ This completes the proof for odd $p.$

Finally, if $p=2$ then our hypothesis ensures $\mathbb{Q}(\zeta_{4n_{2^\prime}})\subseteq F.$ Now
$\mathbb{Q}(4n_{2^\prime})/\mathbb{Q}(\zeta_{4n_{2^\prime}})$ is a cyclic Galois extension of $2$-power degree, and the rest of
the proof goes through as above.
\end{proof}

In the following lemma, we use the notation $\omega_\theta(g)=\frac{h_g\theta(g)}{\theta(1)}$ where $\theta$ is a possibly
reducible character of $G.$ As is well known, $\omega_\theta$ is an algebraic integer provided $\theta$ is induced from an
irreducible character of a subgroup of $G.$

\begin{lem}\label{InductionLemma}
Suppose $\theta=\psi^G$ where $\psi\in\Irr(H)$ for a subgroup $H\subseteq G.$ If $\omega_\psi(h)$ is divisible by
$\abs{Aut^0_H(h)}$ for all $h\in H,$ then $\omega_\theta(g)$ is divisible by $\abs{Aut^0_G(g)}$ for all $g\in G.$
\end{lem}

\begin{proof}
Let $g\in G,$ and let $g_1,...,g_r$ be a set of representatives of the $H$-conjugacy classes in $g^G\cap H.$ Then
$\omega_\theta(g)=\sum_k\omega_{\psi}(g_k).$ (This shows that $\omega_\theta$ is an algebraic integer, as mentioned above.) Let
$n$ be the order of $g,$ and fix a prime $p$ dividing $n.$ We define a field $F$ and extension fields $E_k$ for $1\leq k\leq r$
by
$$F=\mathbb{Q}(\zeta_n)^{Aut^0_{G,p}(g)}\text{ and }E_k=\mathbb{Q}(\zeta_n)^{Aut^0_{H,p}(g_k)}.$$ Note that
$\mathbb{Q}(\zeta_{pn_{p^\prime}},\psi(g_k))\subseteq E_k\subseteq\mathbb{Q}(\zeta_n),$ and if $p=2$ and $4\mid n,$ then also
$i\in E_k.$ Now $Aut^0_{G,p}(g)$ permutes the set of $H$-classes in $g^G\cap H$ and the stabilizer of the class containing $g_k$
is exactly $Aut^0_{H,p}(g_k).$ Hence, if $X$ is a subset of $\set{1,..,r}$ such that the classes $g^H_k$ for $k\in X$ are
representatives for this permutation action, then $$\omega_\theta(g)=\sum_{k\in X}tr_{E_k/F}(\omega_\psi(g_k)).$$ By hypothesis,
$\omega_\psi(g_k)$ is divisible by $\vert Aut^0_{H,p}(g_k)\vert$ for each $k.$ Hence by Lemma \ref{tracelemma},
$tr_{E_k/F}(\omega_\psi(g_k))$ is divisible by $\vert Aut^0_{H,p}(g_k)\vert\abs{E_k:F}=\vert Aut^0_{G,p}(g)\vert.$ The prime $p$
was any prime dividing $n,$ so the result follows.
\end{proof}

\begin{thm}\label{conjwheninduction}
Let $G$ be a finite group and let $p$ be a prime. Let $\chi\in\Irr(G)$ and suppose there exists an expression
$$\chi=\sum_ia_i\lambda_i^G,$$ where $H_i\subseteq G$ are subgroups with $\chi(1)_p$ dividing $\abs{G:H_i}$ and for each $i,$
$\lambda_i$ is a linear character of $H_i$ and $a_i\in\mathbb{Z}.$ Then Conjecture \ref{conj2} holds for $\chi$ at $p.$
\end{thm}

\begin{proof}
Let $g\in G.$ We must prove that $\omega_\chi(g)$ is divisible by $\abs{Aut^0_G(g)}_p.$ This is equivalent to
$$\frac{h_g\chi(g)}{\vert Aut^0_G(g)\vert}=0\mod \chi(1)_p,$$ where the left hand side is an integer because $Aut^0_G(g)$ is a subgroup
of $N_G\langle g\rangle/C_G(g).$ We have
$$\frac{h_g}{\vert Aut^0_G(g)\vert}\chi(g)=\sum_ia_i\frac{h_g}{\vert Aut^0_G(g)\vert}\lambda^G(g).$$ Since $\lambda_i(1)=1$, the conjecture
holds for $\omega_{\lambda_i},$ so by Lemma \ref{InductionLemma}, $\omega_{\lambda_i^G}(g)$ is divisible by $\vert
Aut^0_G(g)\vert,$ or equivalently, $\frac{h_g}{\vert Aut^0_G(g)\vert}\lambda_i^G(g)$ is divisible by $\lambda_i^G(1).$ But
$\lambda_i^G(1)=\abs{G:H_i}$ is divisible by $\chi(1)_p,$ and the result follows.
\end{proof}

\begin{cor}\label{conjsolvable}
Let $G$ be any finite group and let $\chi\in\Irr(G).$ \\$\mathrm{(i)}$ If $G$ is solvable, then Conjecture \ref{conj2} holds for
$\chi.$ \\$\mathrm{(ii)}$ If $G$ is arbitrary but $\chi$ has height zero at a prime $p,$ then Conjecture \ref{conj2} holds for
$\chi$ at $p.$
\end{cor}

\begin{proof}
Immediate by Theorem \ref{conjwheninduction} and Theorem \ref{BrauerInd} when $G$ is solvable, or Theorem \ref{Willems} when
$\chi$ has height zero at $p.$
\end{proof}

Corollary \ref{conjsolvable} provides some evidence for Conjecture \ref{conj2}, and we hope it has some independent interest,
but it does not help us to establish Conjecture \ref{conj1}. Indeed, the hypothesis of Theorem \ref{conjwheninduction} involves
requiring that Conjecture \ref{conj1} holds for $G$ and $\chi$ at $p,$ since (as we already remarked below Theorem
\ref{BrauerInd}), if $\chi$ has an expression of the form $\chi=\sum_ia_i\lambda_i^G$ with $\chi(1)_p$ dividing each
$\abs{G:H_i},$ then $\chi(g)=0$ implies that $g$ is conjugate to an element of some $H_i,$ and then the order of $g$ divides
$\abs{H_i}$ which divides $\abs{G}/\chi(1)_p.$

An additional case where it is easy to prove Conjecture \ref{conj2} is where $\chi$ is the nontrivial irreducible component of a
doubly transitive permutation character:

\begin{thm}\label{perm1}
Suppose $G$ acts doubly transitively on a set $\Omega$ with permutation character $\pi.$ Let $\chi=\pi-1_G,$ so that $\chi$ is
an irreducible character of $G.$ Then Conjecture \ref{conj2} holds for $\chi.$
\end{thm}

\begin{proof}
Fix $g\in G$ and a prime $p.$ Let $Q$ be a Sylow $p$-subgroup of $N_G\langle g_p\rangle\cap C_G(g_{p^\prime}).$ It suffices to
show that $\abs{Q}$ divides $\abs{G}\chi(g)/\chi(1).$ This is equivalent to Conjecture \ref{conj2} except for $p=2,$ where it is
stronger than Conjecture \ref{conj2} as we do not need to take a subgroup of index $2$ in this case. If $H$ is the stabilizer of
a point of $\Omega,$ and $s\in G\backslash H$ is any element then
$$\frac{\abs{G}}{\chi(1)}=\frac{\abs{G}}{\abs{\Omega}-1}=\frac{\abs{G}}{\abs{G:H}-1}=\abs{G:H}\abs{H^s\cap H}.$$

We may assume that $p$ divides $\chi(1)$ or the result is obvious. Hence $\abs{\Omega}$ is prime to $p,$ so $Q$ has a fixed
point $x\in \Omega.$ Let $H$ be the stabilizer of $x,$ so $Q\subseteq H.$ If $g\notin H$ then as $g$ normalizes $Q,$ we find
$Q\subseteq H^g\cap H$ and the result follows from the above formula with $s=g.$ Otherwise $g\in H$ and
$\chi(g)=(\Omega\backslash x)^g.$ Now $Q$ normalizes $\langle g\rangle,$ and so acts on the set $(\Omega\backslash x)^g.$ If
$s\notin H$ then the stabilizer in $Q$ of $y=sx$ is $Q\cap H^s\cap H,$ so $\abs{H^s\cap H}\abs{y^Q}=\abs{H^s\cap
H}\abs{Q}/\abs{Q\cap H^s\cap H}=0\mod\abs{Q}.$ Since $\abs{H^s\cap H}$ is independent of $s$ provided $s\notin H,$ this shows
that $\abs{H^s\cap H}\chi(g)$ is divisible by $\abs{Q},$ and the result follows by the formula for $\abs{G}/\chi(1)$ above.
\end{proof}

Note that by Theorem \ref{1implies2}, Conjecture \ref{conj1} also holds for $\chi$ in this case. Of course this is also easy to
verify directly.

\end{document}